\documentclass[12pt]{amsart}
\usepackage{amsmath}
\usepackage{amssymb}
\usepackage{amsthm}
\usepackage{color}
\usepackage{tikz}
\usepackage[hang,small,bf]{caption}
\usepackage[subrefformat=parens]{subcaption}
\makeatletter
\renewenvironment{proof}[1][\proofname]{\par
  \pushQED{\qed}%
  \normalfont \topsep6\p@\@plus6\p@\relax
  \trivlist
  \item\relax
  {#1\@addpunct{.}}\hspace\labelsep\ignorespaces
}{%
  \popQED\endtrivlist\@endpefalse
}
\makeatother

\newtheorem{theorem}{Theorem}

\theoremstyle{definition}

\theoremstyle{remark}
\newtheorem{definition}{\textup{Definition}}
\newtheorem{remark}[theorem]{\textup{Remark}}

\numberwithin{equation}{section}

\newcommand{\R}{\mathbb R}

\newcommand{\C}{\mathbb C}

\renewcommand{\S}{\mathbb S}

\newcommand{\SO}{{\rm SO}(3)}
\newcommand{\ISO}{{\rm SO}(2)}
\newcommand{\so}{\mathfrak {so}(3)}

\newcommand{\m}{\mathfrak m}
\newcommand{\h}{\mathfrak h}

\newcommand{\del}{\partial}

\newcommand{\g}{\mathfrak{g}}
\newcommand{\la}{\langle}
\newcommand{\ra}{\rangle}

\newcommand{\Isom}{\operatorname{Isom}}

\newcommand{\AS}{\Tilde{\nabla}}

\numberwithin{equation}{section}

\begin{document}
\title[Homogeneous Structures on $\mathbb{S}^2\times \R$ and $\mathbb{H}^2\times \R$]{Homogeneous Structures on $\S^2\times \R$ and $\mathbb{H}^2\times \R$}
%\author{Shimpei Kobayashi}
%\address{Department of Mathematics, Hokkaido University, 
%Sapporo, 060-0810, Japan}
%\email{shimpei@math.sci.hokudai.ac.jp}
%\thanks{The first named author is partially supported by Kakenhi 22K03265.}
%\keywords{Statistical manifolds; conjugate symmetries; $\alpha$-connections; solvable Lie groups; multivariate normal distributions}

%    author two information
\author{Yu Ohno}
\email{ono.yu.h4@elms.hokudai.ac.jp}
\thanks{The author is supported by the establishment of university fellowships towards the creation of
science technology innovation JPMJFS2101.}
%    \subjclass is required.
\dedicatory{}
%    Abstract is required.
\begin{abstract}
 We determine all the homogeneous structure tensors on $\S^2\times \R$ and $\mathbb{H}^2\times \R$. This work together with previous articles \cite{Abe,C,Ca,Seki,TV} yields a complete classification of all the homogeneous structure tensors on three-dimensional homogeneous Riemannian manifolds.

% \medskip
% {\color{magenta}This work together with previous articles 
% \cite{Abe,C,Ca,Seki,TV} yield a complete classification of 
% Grassmann geometry of all the homogeneous structure tensors on three-dimensional homogeneous Riemannian manifolds.}

\end{abstract}

\maketitle
\section*{Introduction}

Locally symmetric Riemannian manifolds are 
characterized by E.~Cartan as Riemannian manifolds whose curvature is 
constant along under parallel translation with respect to the Levi-Civita connection. Ambrose and Singer extended this characterization to homogeneity of Riemannian manifolds in \cite{AS}. 
They proved that a Riemannian manifold $(M,g)$ is locally 
homogeneous if and only if there exists a tensor field $T$ of type $(1,2)$ on $M$ such that
\begin{align}\label{AS}
        \AS g=0,\quad\AS R=0,\quad \AS T=0,
\end{align}
where $\AS$ is an affine connection on $M$ defined by $\AS =\nabla -T$, $\nabla$ is a 
Levi-Civita connection of the metric $g$, and $R$ is a curvature tensor field of $
\nabla$. Such a tensor field $T$ is called a \emph{homogeneous Riemannian structure tensor} \cite{C,TV}. 
Note that Riemannian symmetric spaces have homogeneous Riemannian structure tensor $T=0$. 

% The classification of the homogeneous structures on two-dimensional Riemannian manifold is given by Tricerri and Vanhecke in \cite{TV}. They proved that if $(M,g)$ is a two-dimensional Riemannian manifold and there exists a nonzero homogeneous Riemannian structure $T$, then $(M,g)$ has constant negative curvature. Moreover if $(M,g)$ is a connected, complete and simply connected two-dimensional manifold, then $(M,g)$ admits a homogeneous structure $T\neq 0$ if and only if $(M,g)$ is isometric to the hyperbolic plane $(\mathbb{H}^2,g)$.

Our purpose is to give the classification of the homogeneous structure tensors on 
three-dimensional simply connected Riemannian manifold. In \cite{Seki}, 
it was shown that any three-dimensional connected, simply connected and 
complete homogeneous Riemannian manifold is isometric to either a 
symmetric space; $\S^3$, $\mathbb{H}^3$, $\mathbb{R}^3$, $\S^2\times\R$, $\mathbb{H}^2\times\R$,
or a Riemannian group, \textit{i.e.}, Lie group with left-invariant metric. In \cite{Abe}, Abe determined all the homogeneous structure tensors on 
three-dimensional space forms $\S^3$, $\mathbb{H}^3$ and $\R^3$. For Riemannian Lie groups, 
all the homogeneous structure tensors on Heisenberg group and Berger sphere were determined
in \cite{TV} and \cite{GO}, respectively.
Recently, 
all the homogeneous structures on non-symmetric three dimensional Riemannian Lie groups were determined in \cite{Ca}.

% \medskip

% {\color{magenta}
% In \cite{Abe}, Abe determined all the homogeneous structure tensors on 
% three-dimensional space forms $\S^3$, $\mathbb{H}^3$ and $\R^3$. For Riemannian Lie groups, 
% all the homogeneous structure tensors on Heisenberg group and Berger sphere were determined
% in \cite{TV} and \cite{GO}, respectively.
% Recently, 
% all the homogeneous structures on non-symmetric three dimensional Riemannian Lie groups were determined in \cite{Ca}. 
% }

% \medskip

In this paper, we will determine all the homogeneous structure tensors on 
the Riemannian products $\mathbb{S}^2\times \R$ and $\mathbb{H}^2\times \R$. 
It is easy to see that $\mathbb{S}^2\times \R$ has only one coset space representation 
$\mathbb{S}^2\times \R=\operatorname{SO}(3)\times\R/\operatorname{SO}(2)$ and $\mathbb{H}^2\times \R$ has two coset space representations, 
$\mathbb{H}^2\times \R=\operatorname{SL}(2,\R)\times \R/\operatorname{SO}(2)$ or $\mathbb{H}^2\times \R=\mathbb{H}^2\times \R/\{\mathrm{Id}\}$. 
In the latter case, $\mathbb{H}^2\times \R$ is regarded as a solvable Riemannian group. 
Therefore, we first determine the reductive decomposition corresponding to the above coset space representations, and then we obtain the homogeneous structure tensors from these reductive decompositions. The main results of this paper are as follows.

\begin{theorem}\label{thm:StwoR}
 The Riemannian product $\mathbb{S}^2\times\mathbb{R}$ admits only one type of homogeneous structure tensors up to isomorphism{\rm:}
\begin{align*}
T^{\lambda}=\lambda (dy\otimes dV_{\S^2})\quad (\lambda \ge 0),
\end{align*}
where $dV_{\mathbb{S}^2}$ is the area form of $\mathbb{S}^2$ 
and $y$ is the coordinate of $\mathbb{R}$.
The corresponding coset space representation of each $T^\lambda$ is 
$\mathrm{SO}(3)\times \R/\mathrm{SO}(2)$. 
\end{theorem}

\begin{theorem}\label{thm:HtwoR}
    The Riemannian product $\mathbb{H}^2\times \R$ 
	of the upper half plane model 
	$\mathbb{H}^2$ of the 
	hyperbolic plane with 
	Poincare metric $(dx^2+dy^2)/y^2$ and the real line $(\mathbb{R},dz^2)$ admits two types of 
	homogeneous structure tensors up to isomorphism{\rm:}
    \begin{align*}
    T^{\lambda}&=\lambda \theta^3 \otimes (\theta^1 \wedge \theta^2)=\lambda (\theta_3 \otimes dV_{\mathbb{H}^2}) \quad (\lambda\ge 0), \\
    T&=\theta^1\otimes (\theta^1 \wedge \theta^2)
	=\theta^1\otimes dV_{\mathbb{H}^2},
    \end{align*}
    where $\theta^1=\frac{dx}{y}$, $\theta^2=\frac{dy}{y}$, $\theta^3=dz$, and 
	$dV_{\mathbb{H}^2}$ is the area form of $\mathbb{H}^2$. The corresponding coset space representations 
	of $T^{\lambda}$ and $T$ are $\mathbb{H}^2\times \R=\operatorname{SL}(2,\R)\times \R/\operatorname{SO}(2)$ 
	and $\mathbb{H}^2\times \R=\mathbb{H}^2\times \R/\{\mathrm{Id}\}$, respectively.
\end{theorem}
Here, we denote the covariant tensor field $g(T_XY,Z)$ by the same letter $T$ for simplicity.
Combining Theorem \ref{thm:StwoR} and \ref{thm:HtwoR} with previous studies (\cite{Abe},\cite{C},\cite{Ca},\cite{Seki},\cite{TV}), we can obtain all the homogeneous structure tensors on three-dimensional simply connected Riemannian manifolds. 

\begin{remark}{\rm
Let $M_1$ and $M_2$ be homogeneous Riemannian manifolds, where $M_1$ (resp. $M_2$) has a homogeneous structure tensor $T_1$ (resp. $T_2$). Then, the direct product manifold $M_1\times M_2$ has a homogeneous structure tensor $T=T_1+T_2$. However, since $\S^2$ and $\R$ has only one homogeneous structure tensor $T=0$, the homogeneous structure tensor $T^{\lambda}$ in Theorem \ref{thm:StwoR} cannot be obtained in this way.
}
\end{remark}

\section{Preliminaries}\label{sc:Pre}
In this section, we introduce the construction of homogeneous structure tensor from reductive decomposition.
% Let $(M,g)$ be a connected Riemannian manifold. The Riemannian curvature $R$ is defined by
% \begin{align}\label{def:R}
%     R(X,Y)=\nabla_{[X,Y]}-[\nabla_X,\nabla_Y],
% \end{align}
% where $\nabla$ is the Levi-Civita connection of $g$.
% \begin{definition}
%     A \textit{homogeneous structure} $T$ on $(M,g)$ is a tensor field of type $(1,2)$ which satisfies 
%     \begin{align}
%         \AS g=0,\quad\AS R=0,\quad \AS T=0,
%     \end{align}
%     where $\AS$ is an affine connection on $M$ defined by $\AS =\nabla -T$.
% \end{definition}
% {\color{red}We denote the 
% covariant tensor field
% \[
% g(T(X,Y),Z)
% \]
% by the same letter $T$ for simplicity.
% }

Let $(M,g)$ be a homogeneous Riemannian manifold. Then $M\cong G/H$, where $G$ is 
a connected Lie group acting transitively and effectively on $M$ as a group of isometries and $H$ is the isotropy group at a point $o\in M$, which will be called the \textit{origin} of $M$. We denote by $\g$ and $\h$ the Lie algebras of $G$ and $H$,
 respectively. 
 Since any homogeneous Riemannain manifold is reductive (see p.36 in \cite{C}), 
 there exists a linear subspace $\m$ of $\g$ which is $\operatorname{Ad} (H)$-invariant.
 The linear subspace $\mathfrak{m}$ is called a 
 \emph{Lie subspace} \cite{O}. If $H$ is connected, then the $\operatorname{Ad} (H)$-invariant condition is equivalent to 
$[\m,\h]\subset\m$. For $X\in \g$, we define a vector field $X^*$ on $M$ by
\begin{align}\label{eq:canonical}
    X^*_p=\frac{d}{dt}\exp (tX)\cdot p \biggr|_{t=0}.
\end{align}
Note that $[X^*,Y^*]=-[X,Y]^*$. Then we can identify the Lie subspace $\m$ with the tangent space $T_oM$ via the isomorphism $\tau$ defined by
\begin{align}\label{def:tau}
    \tau (X)=X^*_o \quad (X\in \m).
\end{align}
Then an affine connection $\AS$ of $(M,g)$ corresponding to the reductive decomposition $\g=\m \oplus \h$ is defined by
\begin{align*}
    (\AS_{X^*}Y^*)_o=-([X,Y]_{\m})_o^* \quad (X,Y\in \m).
\end{align*}
One can see that $\AS$ satisfies the equations \eqref{AS}, and then we can obtain a homogeneous structure tensor $T$ by $T=\nabla - \AS$. Moreover, for $u,v,w \in T_o(M)$, we have
\begin{align}
    2g(T_uv,w)_o=g([X,Y]^*,Z^*)_o-g([Y,Z]^*,X^*)+g([Z,X]^*,Y^*) \label{eq:T}
\end{align}
where we take $X, Y, Z \in \m$ such that $\tau (X)=u, \tau (Y)=v,$ and $\tau (Z)=w$. 

% Conversely, let $(M,g)$ be a simply connected Riemannian manifold with a homogeneous Riemannian structure $T$. Fix $o\in M$ and put $\Tilde{\m}= T_oM$. Let $\AS$ be the Ambrose-Singer connection and $\Tilde{R}$ be the curvature tensor of $\AS$. Denote the holonomy algebra of $\AS$ by $\Tilde{\h}$, and then $\Tilde{\h}$ is generated by the operators $\Tilde{R}_{XY}$ for $X,Y \in \Tilde{\m}$. Here, we define a Lie bracket on the direct product $\Tilde{\g}=\Tilde{\m}\oplus \Tilde{\h}$ by
% \begin{align}\label{def:Liebra}
% \begin{split}
%      [U,V]&=UV-VU,\\
%     [U,X]&=U(X),\\
%     [X,Y]&=R(X,Y)+T_XY-T_YX
% \end{split}
% \end{align}
% for all $U,V\in \Tilde{\h}$ and $X,Y\in \Tilde{\m}$. Let $\Tilde{G}$ be the simply connected Lie group whose Lie algebra is $\Tilde{\g}$ and let $\Tilde{H}$ be the connected Lie subgroup of $\Tilde{\h}$ whose Lie algebra is $\Tilde{\h}$. Then $\Tilde{G}$ acts transitively and isometrically on $M$ and $M$ is a coset space $\Tilde{G}/\Tilde{H}$. Let $\Gamma=\{g\in \Tilde{G}|g\cdot x=x \text{~for~all~}x\in M\}$. Then $\Gamma$ is a discrete normal subgroup of $\Tilde{G}$ and $G=\Tilde{G}/\Gamma$ acts transitively, effectively and isometrically on $M$. The isotoropy subgroup $H$ of $G$ at $o \in M$ is $H=\Tilde{H}/\Gamma$. Therefore, we obtain a coset representation $M=G/H$ and $(M,g)$ is a homogeneous Riemannian manifold.

As a special case (\cite{C},\cite{TV}), let $(G,g)$ be a Riemannian group. Then we have a canonical coset representation $G=G/\{e\}$ and reductive decomposition $\g=\{0\}\oplus \g$. Moreover, a homogeneous structure tensor $T$ is determined by the reductive decomposition and $T$ satisfies
\begin{align}
    2g(T_XY,Z)=g([X,Y],Z)-g([Y,Z],X)+g([Z,X],Y) \label{eq:LieT}
\end{align}
for all $X,Y,Z\in \g$.
The connection $\tilde{\nabla}=\nabla-T$ is the so-called 
Cartan-Schouten's $(-)$-connection.

To close this section, we define an isomorphism of homogeneous structure tensors.

\begin{definition}[\cite{TV}]
    Let $T$ be a homogeneous structure tensor on $(M,g)$ and $T'$ be a homogeneous structure tensor on $(M',g')$. Then $T$ and $T'$ are said to be \textit{isomorphic} if there exists an isometry $\phi \colon (M,g) \to (M',g')$ such that
    \begin{align*}
        \phi_*(T_XY)=T'_{\phi_*X}\phi_*Y
    \end{align*}
    for all $X,Y\in \mathfrak{X}(M)$, where $\phi_*$ denotes the differential of $\phi$.
\end{definition}

% A characterization of isomorphic homogeneous structures is given by the following theorem.

% \begin{theorem}\cite{TV}
%     Let $(M,g)$ be a connected, complete and simply connected Riemannian manifolds with homogeneous structures $T$ and $T'$ on $M$. Further let $\g =\h \oplus \m$ (resp. $\g'=\h' \oplus \m'$) be a reductive decomposition determined by $T$ (resp. $T'$). Then the homogeneous Riemannian structures $T$ and $T'$ are isomorphic if and only if there exists a Lie algebra isomorphism $\psi\colon \g \to \g'$ such that $\psi (\h)=\h'$, $\psi (\m)=\m'$ and $\psi |_{\m}\colon \m \to \m'$ is an isometry.
% \end{theorem}

For more details we refer the reader to \cite{TV}.

%{\color{blue}Mention \cite{AS} and \cite{TV} appropriately}

\section{Homogeneous structures on $\mathbb{S}^2\times \R$}

Let us realize the Riemannian product 
$\mathbb{S}^2\times\R$ as a hyperquadric
\[
\mathbb{S}^2\times \R=\{(x^1,x^2,x^3,y)\in \R^4\>|\>(x^1)^2+(x^2)^2+(x^3)^2=1\}
\]
of the Euclidean $4$-space $\R^4$. The metric induced from $\R^4$ 
coincides with the product metric $g=g_{\mathbb{S}^2}+dy^2$.

% \begin{theorem}
% The product manifold 
% {\color{red}$\mathbb{S}^2\times\mathbb{R}$} 
% admits one type of homogeneous structure up to isomorphism{\rm:}
% \begin{align*}
% T^{\lambda}=\lambda (dy\otimes dV_{S^2})\quad (\lambda \ge 0),
% \end{align*}
% {\color{red} where $dV_{\mathbb{S}^2}$ is the area element of $\mathbb{S}^2$. 
% The corresponding coset space representation of each $T^\lambda$ is 
% $(\mathrm{SO}_3\times \R)/\mathrm{SO}_2$. The homogeneous structure 
% is of type $\mathcal{T}_2\oplus\mathcal{T}_3$. 
% In particular $T^{\lambda}$ is of type $\mathcal{T}_3$ 
% if and only if $\lambda=0$. For $\lambda\not=0$, 
% $T^{\lambda}$ is not of type $\mathcal{T}_2$.
% }
% \end{theorem}
\begin{proof}[Proof of Theorem 1]
The identity component $\Isom _0(\mathbb{S}^2\times \R)$ 
of the full isometry group of $\mathbb{S}^2\times\R$ is 
$\Isom _0(\mathbb{S}^2\times \R)=\SO \times \R$ with Lie algebra is 
$\so \oplus \R$. Under the hyperquadric model of $\mathbb{S}^2\times\R$, 
the isometry group $\SO \times \R$ is identified with
\[
\left\{
\left.
\left(
\begin{array}{cc}
A & 0\\
0 & e^{t}
\end{array}
\right)
\>\right|
\>a\in\SO,\>
t\in\mathbb{R}
\right\}\subset\mathrm{GL}(4,\R).
\]

We denote the isotropy subgroup of $\SO \times \R$ at 
the origin $o=(1,0,0,0) \in \mathbb{S}^2 \times \R$ by 
$H$, which is isomorphic to $\ISO$. 
Since there are no subgroups of dimension $2$ in $\SO$, 
the connected subgroup $G$ of $\Isom _0
(\mathbb{S}^2 \times \R)$ acting transitively on 
$\mathbb{S}^2 \times \R$ is $\SO \times \R$. 
We take an orthonormal basis $\{u_1, u_2, u_3, e\}$ of $\so \oplus \R$ given by 
\begin{align}
    u_1=E_{23}-E_{32},\quad u_2=E_{12}-E_{21},\quad u_3=E_{31}-E_{13},\quad e=E_{44}, \label{def:u} 
\end{align}
where $E_{ij}$ is the matrix unit of $4\times 4$ matrix. Then, we have
\begin{align}
    [u_1,u_2]=u_3,\quad [u_2,u_3]=u_1,\quad [u_3,u_1]=u_2, \quad [u_i,e]=0\quad (i=1,2,3). \label{bra:u}
\end{align}
Here, the Lie algebra of $H$ is $\mathfrak {h}=\la\{u_1\}\ra$. Moreover, from \eqref{eq:canonical}, \eqref{def:tau} and \eqref{def:u}, we have
\begin{align}
    \tau(u_2)=\left.\frac{\del}{\del x^2}\right| _{o},
	\quad 
	\tau(u_3)=-\left.\frac{\del}{\del x^3}\right| _{o},\quad 
	\tau(e)=\left.\frac{\del}{\del y}\right| _{o}. \label{cal:u}
\end{align}

Let $\mathfrak{m}$ be a Lie subspace
of $\mathfrak{h}$ in $\so \oplus \R$, 
then $\mathfrak{m}$ is expressed as $\la \{ u_2+\lambda_2u_1, u_3+\lambda_3u_1, e+\lambda u_1\} \ra$ for $\lambda_2,\lambda_3,\lambda\in \R$. 
Moreover, from $[\mathfrak{m},\mathfrak{h}_p]\subset \mathfrak{m}$, it follows that $\lambda_2=\lambda_3=0$, and therefore $\mathfrak{m}$ is given by
\begin{align*}
    \mathfrak{m}=\m ^{\lambda}=\la \{ u_2,u_3,e+\lambda u_1\} \ra \quad (\lambda\in \R). 
\end{align*}
Then, from \eqref{eq:T}, \eqref{bra:u} and \eqref{cal:u}, the homogeneous structure $T^{\lambda}$ defined from the reductive decomposition $\so \oplus \R=\m^{\lambda} \oplus \mathfrak{h}$ satisfies

    \begin{align*}
        \left\{
        \begin{array}{llll}
            \left.T^{\lambda}_{\frac{\partial}{\partial x^2}}\frac{\partial}{\partial x^2}\right|_{o}=0, 
			& \left.T^{\lambda}_{\frac{\partial}{\partial x^2}}\frac{\partial}{\partial x^3}\right|_{o}=0, 
			&\left.T^{\lambda}_{\frac{\partial}{\partial x^2}}\frac{\partial}{\partial y}\right|_{o}=0, \\
            \left.T^{\lambda}_{\frac{\partial}{\partial x^3}}\frac{\partial}{\partial x^2}\right|_{o}=0, & \left.T^{\lambda}_{\frac{\partial}{\partial x^3}}\frac{\partial}{\partial x^3}\right|_{o}=0, 
			&\left.T^{\lambda}_{\frac{\partial}{\partial x^3}}\frac{\partial}{\partial y}\right|_{o}=0, \\
            \left.T^{\lambda}_{\frac{\partial}{\partial y}}\frac{\partial}{\partial x^2}\right|_{o}
			=\left.\lambda\frac{\partial}{\partial x^3}\right|_o, &\left.T^{\lambda}_{\frac{\partial}{\partial y}}\frac{\partial}{\partial x^3}\right|_{o}
			=-\left.\lambda\frac{\partial}{\partial x^2}\right|_{o}, 
			&\left.T^{\lambda}_{\frac{\partial}{\partial y}}\frac{\partial}{\partial y}\right|_{o}=0. 
        \end{array}
        \right.
    \end{align*}

Therefore, we obtain
\begin{align*}
    T^{\lambda}|_{o}=\lambda (dy\otimes (dx^2\wedge dx^3))|_{o}
	=\lambda (dy  \otimes dV_{\mathbb{S}^2})|_{o},
\end{align*}
where $dV_{\mathbb{S}^2}$
denotes the area form of $\mathbb{S}^2$. 
Since $T^{\lambda}$ and $dy\otimes dV_{\mathbb{S}^2}$ are $(\SO \times \R)$-invariant, we have
\begin{align*}
    T^{\lambda}=\lambda (dy\otimes dV_{\mathbb{S}^2}).
\end{align*}
In addition, for an isometry $\phi$, it follows that $\phi^*T^{\lambda}=T^{\lambda}$
if $\phi$ preserves the orientation, and
$\phi^*T^{\lambda}=-T^{\lambda}$
if $\phi$ changes the orientation. Therefore, $T^{\lambda}$ and $T^{\mu}$ are isomorphic if and only if $\lambda=\pm \mu$. This completes the proof.    
\end{proof}

% {\color{magenta}
% ************************************************************************************
% }

% We regard $\{u_2,u_3,e\}$ as an 
% $(\mathrm{SO}_3\times\mathbb{R})$-invariant 
% orthonormal frame field on $\mathbb{S}^2\times\mathbb{R}$. 
% Denote by $\{\theta^2,\theta^2,dy\}$ the 
% dual coframe field. Then 
% \[
% T=\lambda (dy\otimes \theta^2\wedge\theta^3)
% \]for some constant $\lambda$. One can see that 
% $dV_{\mathbb{S}^2}=\theta^2\wedge \theta^3$.

% Let us consider the homogeneous space $\mathbb{S}^2
% =\mathrm{SO}_3/\mathrm{SO}_2$. 
% The tangent space of $\mathbb{S}^2$ at the origin 
% is identified with 
% \[
% \{c_1u_2+c_2u_3\>|\>c_2,c_3\in\mathbb{R}\}.
% \]
% The standard complex structure of $\mathbb{S}^2$ is defined by
% \[
% Ju_2=u_3,\quad Ju_3=-u_2.
% \]
% We extend $J$ to $\mathfrak{m}^{0}$ by the rule $Je=0$. Then 
% one can check that 
% \begin{align*}
% g(T^{\lambda}(X,Y),Z)=&\lambda\,dy(X)(\theta^2\wedge\theta^3)(Y,Z)
% =\frac{\lambda}{2}dy(X)
% \left(g(Y,u_2)g(Z,u_3)-g(Y,u_3)g(Z,u_2)\right)
% \\
% =&\frac{\lambda}{2}dy(X)\,
% g(g(Y,u_2)u_3-g(Y,u_3)u_3,Z)
% =\frac{\lambda}{2}dy(X)\,g(JY,Z).
% \end{align*}
% Thus we obtain 
% \[
% T^{\lambda}(X,Y)=\lambda\,dy(X)JY.
% \]

% We need to comment about the homogeneous structures of $\mathbb{S}^2$.

% {\color{magenta}
% ************************************************************************************
% }
\section{Homogeneous structures on $\mathbb{H}^2\times \R$} 

Let us realize the Riemannian product $\mathbb{H}^2\times \R$ as 
\begin{align*}
    \mathbb{H}^2\times \R=\{(w=x+yi,z)\in \C \times \R \>|\> y>0\}
\end{align*}
with the product metric 
\[g=g_{\mathbb{H}^2}+dz^2=\frac{1}{y^2}(dx^2+dy^2)+dz^2.\]

% \begin{theorem}
%     The product manifold $\mathbb{H}^2\times \R$ admits two type of homogeneous structures up to isomorphism:
%     \begin{align*}
%     T^{\lambda}&=\lambda \theta^3 \otimes (\theta^1 \wedge \theta^2) \quad (\lambda\ge 0), \\
%     T&=\theta^1\otimes (\theta^1 \wedge \theta^2)
% 	{\color{red}=\theta^1\otimes dV_{\mathbb{H}^2},}
%     \end{align*}
%     where $\theta^1=\frac{dx}{y}$, $\theta^2=\frac{dy}{y}$, $\theta^3=dz$, and 
% 	$dV_{\mathbb{H}^2}$ is the area form of $\mathbb{H}^2$.
	
% {\color{red} The corresponding coset space representations are
% \begin{itemize}
% \item $T^{\lambda}${\rm:} $(\mathrm{SL}_2\mathbb{R}\times \R)/ \mathrm{SO}_2$.
% \item $T${\rm:} $(\mathbb{H}^2\times \R)/\{\mathrm{Id}\}$.
% \end{itemize}
% }	
% The homogeneous structure $T^{\lambda}$ is of type $\mathcal{T}_2\oplus \mathcal{T}_3$ and $T$ is of type $\mathcal{T}_1\oplus \mathcal{T}_2$. In particular $T^{\lambda}$ is of type $\mathcal{T}_3$ if and only if $\lambda=0$.
% \end{theorem}

\begin{proof}[Proof of Theorem 2]
   The identity component $\operatorname{Isom}_0(\mathbb{H}^2\times \R)$ of the full isometry group of 
   $\mathbb{H}^2\times\R$ is $\operatorname{Isom}_0(\mathbb{H}^2\times \R)=\operatorname{SL}(2,\R)\times \R$ with Lie algebra 
   $\mathfrak{sl}(2, \R) \oplus \R$. The isometric action of $\operatorname{SL}(2,\R)\times \R$ on $\mathbb{H}^2\times \R$ is defined by
   \begin{align*}
       \left( 
       \begin{pmatrix}
           a & b \\
           c & d
       \end{pmatrix}
       , s
       \right) \cdot
       (w,z)
       =\left(
        \frac{aw+b}{cw+d},z+s
       \right).
   \end{align*}
   Then the isotoropy subgroup $H$ of $\operatorname{SL}(2,\R)\times \R$ at 
   the origin ${o}=(i,0)$ is $\operatorname{SO}(2)$. 
   Let $G$ be a connected subgroup of $\operatorname{SL}(2,\R)\times \R$ acting transitively on $\mathbb{H}^2\times \R$, and then $\dim G=4$ or $3$. 
   
   First, we consider the case $\dim G=4$. This means that $G=\operatorname{SL}(2,\R)\times \R$. 
   We identify the Lie group $\operatorname{SL}(2,\R)\times \R$ with the linear Lie group
   \begin{align*}
       \operatorname{SL}(2,\R)\times \R=\left\{
       \left.
       \begin{pmatrix}
           A & 0\\
           0 & e^t
       \end{pmatrix}
       \right|
       A\in \operatorname{SL}(2,\R), t\in \R
       \right\}
       \subset\operatorname{GL}(3,\R).
   \end{align*}
   Moreover, we take an orthonormal basis $\{v_1,v_2,v_3,e\}$ of 
   $\mathfrak{sl}( 2, \R) \oplus \R$ given by
   \begin{align}
       v_1=\frac{1}{2}(E_{21}-E_{12}),\quad v_2=\frac{1}{2}(E_{12}+E_{21}),\quad v_3=\frac{1}{2}(E_{22}-E_{11}),\quad e=E_{33} \label{def:v}
   \end{align}
   where $E_{ij}$ is the matrix unit of $3\times 3$ matrix. Then, we have
   \begin{align*}
       [v_1,v_2]=v_3,\quad [v_2,v_3]=-v_1,\quad [v_3,v_1]=v_2,\quad [v_i,e]=0\quad (i=1,2,3).
   \end{align*}
   Here, the Lie algebra of $H$ is 
   $\mathfrak{h}=\la\{v_1\}\ra$. From \eqref{eq:canonical}, \eqref{def:tau} and \eqref{def:v}, we have
   \begin{align*}
       \tau(v_2)=\left.\frac{\partial}{\partial x}\right|_p,\quad \tau(v_3)=-\left.\frac{\partial}{\partial y}\right|_p,\quad \tau(e)=\left.\frac{\partial}{\partial z}\right|_p.
   \end{align*}
   Let $\m$ be a Lie subspace of 
   $\mathfrak{h}$ in 
   $\mathfrak{sl}(2,\R)\oplus \R$. 
   As in the case of $\S^2\times \R$, from 
  $[\mathfrak{h},\m]\subset \m$, the Lie subspace $\m$ is given by
   \begin{align*}
       \m=\m^{\lambda}=\la\{v_2,v_3,e+\lambda v_1\}\ra \quad (\lambda \in \R),
   \end{align*}
   and the homogeneous structure tensor $T^{\lambda}$ determined by the reductive decomposition 
   $\mathfrak{sl}(2,\R)\oplus \R=\m^{\lambda} \oplus \mathfrak{h}$ satisfies
   \begin{align*}
       T^{\lambda}=\lambda(dz\otimes dV_{\mathbb{H}^2}).
   \end{align*}
   For the same reason as in $\S^2\times \R$, $T^{\lambda}$ and $T^{\mu}$ are isomorphic if and only if $\lambda=\pm \mu$.

   Next, we consider the case $\dim G=3$. In this case, the coset space representation is $\mathbb{H}^2\times \R=G/\{e\}$, and then $\mathbb{H}^2\times \R$ is identified with a Lie group $G$. Thus, $G$ is isomorphic to the solvable Lie group 
   \begin{align*}
        \left.
        \mathbb{H}^2\times \R=\left\{
        \begin{pmatrix}
        y&x&0\\
        0&1&0\\
        0&0&e^z
        \end{pmatrix}
        \right|
        x,y,z\in\R, y>0
        \right\},
    \end{align*}
   with the left invariant metric
    \begin{align*}
        g=\frac{1}{y^2}(dx^2+dy^2)+dz^2.
    \end{align*}
    We define a left invariant orthonormal vector field $\{e_1,e_2,e_3\}$ as
    \begin{align}\label{def:e}
        e_1=y\frac{\del}{\del x}, \quad e_2=y\frac{\del}{\del y}, \quad e_3=\frac{\del}{\del z}.
    \end{align}    
    Then we have
    \begin{align}\label{bra:e}
        [e_1,e_2]=-e_1,\quad [e_1,e_3]=[e_2,e_3]=0.
    \end{align}
    Moreover, from \eqref{eq:LieT} and \eqref{bra:e}, the homogeneous structure tensor $T$ determined by the reductive decomposition $\g=\{0\}\oplus\g$ satisfies
    \begin{align}\label{cal:T}
        \left\{
        \begin{array}{llll}
            T_{e_1}e_1=e_2, & T_{e_1}e_2=-e_1, &T_{e_1}e_3=0, \\
            T_{e_2}e_1=0, & T_{e_2}e_2=0, &T_{e_2}e_3=0, \\
            T_{e_3}e_1=0, &T_{e_3}e_2=0, &T_{e_3}e_3=0. 
        \end{array}
        \right.
    \end{align}
    Thus, from \eqref{def:e} and \eqref{cal:T}, we obtain
    \begin{align*}
        T=\theta^1\otimes (\theta^1\wedge \theta^2),
    \end{align*}
    where $\theta^1=\frac{dx}{y}$, $\theta^2=\frac{dy}{y}$. This completes the proof.
\end{proof}

\subsection*{Acknowledgements}
We would like to thank Prof. Jun-ichi Inoguchi and Prof. Shimpei Kobayashi comment on the manuscripts and letting us know several related references.


\begin{thebibliography}{99}

\bibitem{Abe}
K.~Abe,
\newblock The classification of homogeneous structures on 3-dimensional space forms,
{\em Math. J. Okayama Univ.} \textbf{28} (1986), 173--189.

\bibitem{AS}
W.~Ambrose and I.M.~Singer,
\newblock On homogeneous Riemannian manifolds,
{\em Duke Math. J.} \textbf{25} (1958), 647--669.

\bibitem{C}
G.~Calvaruso, M.~Castrill{\'o}n L{\'o}pez,
\textit{Pseudo-Riemannian Homogeneous Structures}, 
Developments in Mathematics
vol. 59. Springer, Cham, 2019.

\bibitem{Ca}
E.~Calvi{\~n}o-Louzao, M.~Ferreiro-Subrido, E.~Garc{\'\i}a-R{\'\i}o and R. V{\'a}zquez-Lorenzo,
\newblock Homogeneous Riemannian structures in dimension three,
{\em Rev. Real Acad. Cienc. Exactas Fis. Nat. Ser. A-Mat.} \textbf{117} (2023), 70.

\bibitem{GO}
P.~M.~Gadea and J.~A.~Oubi{\~n}a,
\newblock Homogeneous Riemannian structures on Berger 3-spheres,
{\em Proc. Edinb. Math. Soc.} \textbf{48} (2005), no.~2, 375--387.

\bibitem{O}
B.~O'Neill,
\textit{Semi-Riemannian Geometry with Applications to Relativity},
Academic press, 1983.

\bibitem{Seki}
K.~Sekigawa,
\newblock On some $3$-dimensional curvature homogeneous spaces,
{\em Tensor N.S.} \textbf{31} (1977), 87--97.

\bibitem{TV}
F.~Tricerri and L.~Vanhecke,
\textit{Homogeneous structures on Riemannian manifolds},
Cambridge university press, 1983.



\end{thebibliography}
\end{document}